\documentclass[11pt]{article}

\usepackage{amsmath}
\usepackage{amsthm}
\usepackage{amssymb}
\usepackage{latexsym}
\usepackage{pstricks}
\usepackage{graphicx, pst-plot, pst-node, pst-text, pst-tree}
\usepackage{titlefoot}
\usepackage{enumerate}
\usepackage{titlesec}
\usepackage[small,it]{caption}
\usepackage{color}
\definecolor{lightgray}{rgb}{0.9, 0.9, 0.9}
\definecolor{darkgray}{rgb}{0.7, 0.7, 0.7}
\definecolor{darkblue}{rgb}{0, 0, .4}

\DeclareFontFamily{U}{noplusfont}{}
\DeclareFontShape{U}{noplusfont}{m}{n}
   {  <10><12> noplusfont10}{}
\newcommand{\noplus}{\textrm{\usefont{U}{noplusfont}{m}{n}\symbol{'010}}}
\newcommand{\nominus}{\textrm{\usefont{U}{noplusfont}{m}{n}\symbol{'011}}}

\usepackage[bookmarks,breaklinks]{hyperref}
\hypersetup{
	colorlinks=true,
	linkcolor=darkblue,
	anchorcolor=darkblue,
	citecolor=darkblue,
	pagecolor=darkblue,
	urlcolor=darkblue,
	pdfpagemode=UseThumbs,
	pdftitle={Simple permutations and algebraic generating functions},
	pdfsubject={Combinatorics},
	pdfauthor={Robert Brignall, Sophie Huczynska, and Vincent Vatter},
	pdfkeywords={permutation class, restricted permutation, simple permutation}
}

\newtheorem{theorem}{Theorem}[section]
\newtheorem{proposition}[theorem]{Proposition}
\newtheorem{lemma}[theorem]{Lemma}

\newtheorem{question}[theorem]{Question}

   \newtheoremstyle{example}{\topsep}{\topsep}%
     {}
     {}
     {\bfseries}
     {.}
     {.5em}
     {\thmname{#1}\thmnumber{ #2}:\thmnote{ #3}}

   \theoremstyle{example}
   \newtheorem{example}[theorem]{Example}

\newcounter{todocounter}

\newcommand{\minisec}[1]{\bigskip\noindent{\bf #1.}}

\makeatletter
\newfont{\footsc}{cmcsc10 at 8truept}
\newfont{\footbf}{cmbx10 at 8truept}
\newfont{\footrm}{cmr10 at 10truept}
\renewenvironment{abstract}%
		{
		  \begin{list}{}%
		     {\setlength{\rightmargin}{1in}%
		      \setlength{\leftmargin}{1in}}%
		   \item[]\ignorespaces\begin{small}}%
		 {\end{small}\unskip\end{list}}

\datefoot{\today}
\amssubj{05A15, 05A05}

\newpagestyle{main}[\small]{
	\headrule
	\sethead[\usepage][][]
	{\sc Simple Permutations and Algebraic Generating Functions}{}{\usepage}}

\title{\sc{Simple Permutations and Algebraic Generating Functions}}
\author{\sc{Robert Brignall, Sophie Huczynska\thanks{Supported by a Royal Society Dorothy Hodgkin Research Fellowship.}, and Vincent Vatter\thanks{Supported by EPSRC grant GR/S53503/01.}}\\
\small School of Mathematics and Statistics\\[-3pt]
\small University of St Andrews\\[-3pt]
\small St Andrews, Fife, Scotland\\[-3pt]
\small \texttt{\{robertb, sophieh, vince\}@mcs.st-and.ac.uk}\\[-3pt]
\small \texttt{http://turnbull.mcs.st-and.ac.uk/\~{}\{\href{http://turnbull.mcs.st-and.ac.uk/~robertb}{robertb}, \href{http://turnbull.mcs.st-and.ac.uk/~sophieh}{sophieh}, \href{http://turnbull.mcs.st-and.ac.uk/~vince}{vince}\}}\\[-10pt]}

\titleformat{\section}
	{\large\sc}
	{\thesection.}{1em}{}	

\date{}

\begin{document}
\maketitle

\pagestyle{main}

\newcommand{\Av}{\operatorname{Av}}
\newcommand{\C}{\mathcal{C}}
\newcommand{\W}{\mathcal{W}}
\newcommand{\Si}{\operatorname{Si}}
\newcommand{\conv}{\operatorname{conv}}
\newcommand{\strongcomp}{\operatorname{sc}}
\newcommand{\Dash}{{\mbox{\small-}}}
\newcommand{\Cong}{\equiv}
\newcommand{\Mod}{\mathop{\rm mod}\nolimits}

\begin{abstract}
A simple permutation is one that does not map a nontrivial interval onto an interval.  It was recently proved by Albert and Atkinson that a permutation class with only finitely simple permutations has an algebraic generating function.  We extend this result to enumerate permutations in such a class satisfying additional properties, e.g., the even permutations, the involutions, the permutations avoiding generalised permutations, and so on.
\end{abstract}

\section{Introduction}\label{sp1-intro}

Substitution decompositions (known also as modular decompositions, disjunctive decompositions, and $X$-joins) have proved to be a useful technique in a wide range of settings, ranging from game theory to combinatorial optimization, see M\"ohring~\cite{mohring:algorithmic-asp:a} or M\"ohring and Radermacher~\cite{mohring:substitution-de:} for extensive references.  Although substitution decompositions are most often applied to algorithmic problems, here we apply them to the enumeration of permutation classes, a.k.a. restricted permutations.

The permutation $\pi$ is said to {\it contain\/} the permutation $\sigma$, written $\sigma\le\pi$, if $\pi$ has a subsequence that is order isomorphic to $\sigma$.  For example, $\pi=491867532$ contains $\sigma=51342$, as can be seen by considering the subsequence $91672$ ($=\pi(2),\pi(3),\pi(5),\pi(6),\pi(9)$), and such a subsequence is called a {\it copy\/} of $\sigma$ in $\pi$.  This pattern-containment relation is a partial order on permutations.  We refer to downsets of permutations under this order as {\it permutation classes\/}.  In other words, if $\C$ is a permutation class, $\pi\in\C$, and $\sigma\le\pi$, then $\sigma\in\C$.  Permutation classes arise naturally in a variety of disparate fields, ranging from the analysis of sorting machines (dating back to Knuth~\cite{knuth:the-art-of-comp:}, who proved that a permutation is stack-sortable if and only if it lies in the class $\Av(231)$) to the study of Schubert varieties (see, e.g., Lakshmibai and Sandhya~\cite{lakshmibai:criterion-for-s:}).  The most frequently investigated property of permutation classes is their enumeration.

We shall denote by $\C_n$ the set $\C \cap S_n$, i.e.\ the permutations in $\C$ of length
$n$, and we refer to $\sum |\C_n| x^n$ as the {\it generating function for $\C$\/}.
Recall that
an {\it antichain} is a set of pairwise incomparable elements.
For any permutation class $\C$, there is a unique (possibly infinite) antichain $B$ such that $\C=\Av(B)=\{\pi: \beta \not \leq\pi\mbox{ for all } \beta \in B\}$. This antichain $B$, which consists of the minimal permutations not in $\C$, is called the {\it basis} of $\C$.

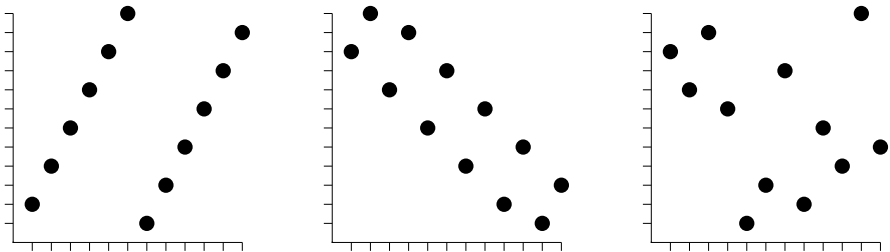
\begin{figure}
\begin{center}
\begin{tabular}{ccccc}
\psset{xunit=0.01in, yunit=0.01in}
\psset{linewidth=0.005in}
\begin{pspicture}(0,0)(120,120)
\psaxes[dy=10,Dy=1,dx=10,Dx=1,tickstyle=bottom,showorigin=false,labels=none](0,0)(120,120)
\pscircle*(10,20){0.04in}
\pscircle*(20,40){0.04in}
\pscircle*(30,60){0.04in}
\pscircle*(40,80){0.04in}
\pscircle*(50,100){0.04in}
\pscircle*(60,120){0.04in}
\pscircle*(70,10){0.04in}
\pscircle*(80,30){0.04in}
\pscircle*(90,50){0.04in}
\pscircle*(100,70){0.04in}
\pscircle*(110,90){0.04in}
\pscircle*(120,110){0.04in}
\end{pspicture}
&\rule{10pt}{0pt}&
\psset{xunit=0.01in, yunit=0.01in}
\psset{linewidth=0.005in}
\begin{pspicture}(0,0)(120,120)
\psaxes[dy=10,Dy=1,dx=10,Dx=1,tickstyle=bottom,showorigin=false,labels=none](0,0)(120,120)
\pscircle*(10,100){0.04in}
\pscircle*(20,120){0.04in}
\pscircle*(30,80){0.04in}
\pscircle*(40,110){0.04in}
\pscircle*(50,60){0.04in}
\pscircle*(60,90){0.04in}
\pscircle*(70,40){0.04in}
\pscircle*(80,70){0.04in}
\pscircle*(90,20){0.04in}
\pscircle*(100,50){0.04in}
\pscircle*(110,10){0.04in}
\pscircle*(120,30){0.04in}
\end{pspicture}
&\rule{10pt}{0pt}&
\psset{xunit=0.01in, yunit=0.01in}
\psset{linewidth=0.005in}
\begin{pspicture}(0,0)(120,120)
\psaxes[dy=10,Dy=1,dx=10,Dx=1,tickstyle=bottom,showorigin=false,labels=none](0,0)(120,120)
\pscircle*(10,100){0.04in}
\pscircle*(20,80){0.04in}
\pscircle*(30,110){0.04in}
\pscircle*(40,70){0.04in}
\pscircle*(50,10){0.04in}
\pscircle*(60,30){0.04in}
\pscircle*(70,90){0.04in}
\pscircle*(80,20){0.04in}
\pscircle*(90,60){0.04in}
\pscircle*(100,40){0.04in}
\pscircle*(110,120){0.04in}
\pscircle*(120,50){0.04in}
\end{pspicture}
\end{tabular}
\end{center}
\caption{The plots of three simple permutations of length $12$.}\label{simple-exs}
\end{figure}

An {\it interval\/} in the permutation $\pi$ is a set of contiguous indices $I=[a,b]$ such that the set of values $\pi(I)=\{\pi(i) : i\in I\}$ also forms an interval of natural numbers.  Every permutation $\pi$ of $[n]=\{1,2,\dots,n\}$ has intervals of length $0$, $1$, and $n$; $\pi$ is said to be {\it simple\/} if it has no other intervals (such intervals are called {\it proper\/}).  Figure~\ref{simple-exs} shows three simple permutations.  Albert and Atkinson~\cite{albert:simple-permutat:} were the first to establish the link between simple permutations and the enumeration of permutation classes; they proved that every permutation class with only finitely many simple permutations has an algebraic generating function.  Our main theorem implies the following result.  (The terms in Theorem~\ref{sp1-main} are defined at the end of this section.)

\begin{theorem}\label{sp1-main}
In a permutation class $\C$ with only finitely many simple permutations, the following sequences have  algebraic generating functions:
\begin{itemize}
\item the number of alternating permutations in $\C_n$,
\item the number of even permutations in $\C_n$,
\item the number of Dumont permutations of the first kind in $\C_n$,
\item the number of permutations in $\C_n$ avoiding any finite set of blocked or barred permutations, and
\item the number of involutions in $\C_n$.
\end{itemize}
Moreover, these conditions can be combined in any finite manner desired.
\end{theorem}

One class to which this theorem applies is $\Av(132)$.  In any permutation from $\Av(132)$, all entries to the left of the maximum must be greater than all entries to the right.  This shows that $\Av(132)$ has only three simple permutations ($1$, $12$, and $21$).  Thus Theorem~\ref{sp1-main} helps to explain why there are known formulas and recurrences for enumerating:
\begin{itemize}
\item $\Av(132,\beta)$ --- Mansour and Vainshtein~\cite{mansour:restricted-132-:a},
\item involutions in $\Av(132,\beta)$ --- Guibert and Mansour~\cite{guibert:restricted-132-:},
\item even involutions in $\Av(132,\beta)$ --- Guibert and Mansour~\cite{guibert:some-statistics:},
\item involutions in $\Av(231,\beta)$ --- Egge and Mansour~\cite{egge:231-avoiding-in:},
\item even permutations in $\Av(132,\beta)$ --- Mansour~\cite{mansour:restricted-even:},
\item permutations in $\Av(132)$ avoiding a blocked permutation --- Mansour~\cite{mansour:restricted-1-3-:},
\item alternating permutations in $\Av(132)$ avoiding a blocked permutation --- Mansour~\cite{mansour:restricted-132-:},
\item Dumont permutations of the first kind in $\Av(132,\beta)$ --- Mansour~\cite{mansour:restricted-132-:b},
\item permutations avoiding $132$, $1\Dash23$, and the permutation or blocked permutation $\beta$ --- Elizalde and Mansour~\cite{elizalde:restricted-motz:},
\item permutations in $\Av(132,\beta)$ that are West-$2$-stack-sortable --- Egge and Mansour~\cite{egge:132-avoiding-tw:}, and
%
%
\end{itemize}
Moreover, in Brignall, Huczynska, and Vatter~\cite{brignall:simple-permutat:a} a decomposition theorem for simple permutations is proved which implies that for any fixed $r$, the class consisting of all permutations with at most $r$ copies of $132$ contains only finitely many simple permutations.  Therefore, Theorem~\ref{sp1-main} also implies that the following sets have algebraic generating functions:
\begin{itemize}
\item permutations with at most $r$ copies of $132$ --- B\'ona~\cite{bona:the-number-of-p:} and Mansour and Vainshtein~\cite{mansour:counting-occurr:},
\item even permutations with at most $r$ copies of $132$ --- Mansour~\cite{mansour:counting-occurr:a},
\item involutions with at most $r$ copies of $231$ --- Mansour, Yan, and Yang~\cite{mansour:counting-occurr:b},
\item alternating permutations with at most $r$ copies of $132$ --- Mansour~\cite{mansour:restricted-132-:}, and
\end{itemize}

We conclude this section by defining the terms involved in Theorem~\ref{sp1-main}.

The permutation $\pi\in S_n$ is {\it alternating\/} if for all $i\in[2,n-1]$, $\pi(i)$ does not lie between $\pi(i-1)$ and $\pi(i+1)$.  A permutation is {\it Dumont of the first kind\/} if each even entry is immediately followed by a smaller entry and each odd entry is either immediately followed by a larger entry or occurs last (this dates back to Dumont~\cite{dumont:interpretations:}).

A {\it barred permutation\/} is a permutation in which one or more of the entries is barred; for $\pi$ to avoid the barred permutation $\sigma$ means that every set of entries of $\pi$ that is order isomorphic to the nonbarred entries of $\sigma$ can be extended to a set order isomorphic to $\sigma$ itself.  For example, $51342$ avoids $\overline{3}12$ because every noninversion (i.e., copy of $12$) can be extended to a copy of $312$ (e.g., prefix the $5$), but $51342$ contains $2\overline{1}3$ because the $1$ and $3$ are order isomorphic to $23$, but there is no way to extend this to a copy of $213$.  Barred permutations have arisen several times in the permutation pattern literature.  For example, under West's notion of $2$-stack sorting (see West~\cite{west:sorting-twice-t:}) the permutations that can be sorted are those that avoid $2341$ and $3\overline{5}241$, while Bousquet-M\'elou and Butler~\cite{bousquet-melou:forest-link-per:} characterise the permutations corresponding to locally factorial Schubert varieties in terms of barred permutations.

A {\it blocked permutation\/} is a permutation containing dashes indicating the entries that need not occur consecutively (in the normal pattern-containment order, no entries need occur consecutively).  For example, $51342$ contains two copies of $3\Dash12$: $513$ and $534$, but note that $514$ is not a copy of $3\Dash 12$ because the $1$ and $4$ are not adjacent.  Babson and Steingr{\'{\i}}msson~\cite{babson:generalized-per:} introduced blocked permutations and showed that they could be used to express most Mahonian statistics.  For example, the major index of $\pi$ is equal to the total number of copies of $1\Dash32$, $2\Dash31$, $3\Dash21$, and $21$ in $\pi$.

\section{Inflations and Wreath Closures}\label{sp1-inflations}

Here we briefly review some required terminology and results.  Given $\sigma\in S_m$ and nonempty permutations $\alpha_1,\dots,\alpha_m$, the {\it inflation\/} of $\sigma$ by $\alpha_1,\dots,\alpha_m$ --- denoted $\sigma[\alpha_1,\dots,\alpha_m]$ --- is the permutation obtained by replacing each entry $\sigma(i)$ by an interval that is order isomorphic to $\alpha_i$.  For example, $2413[1,132,321,12]=479832156$ (see Figure~\ref{fig-479832156}).  Simple permutations cannot be deflated.  Conversely:

\begin{figure}
\begin{center}
\psset{xunit=0.01in, yunit=0.01in}
\psset{linewidth=0.005in}
\begin{pspicture}(0,0)(93,93)
\psaxes[dy=10, Dy=1, dx=10, Dx=1, tickstyle=bottom, showorigin=false, labels=none](0,0)(90,90)
\psframe[linecolor=darkgray,fillstyle=solid,fillcolor=lightgray,linewidth=0.02in](7,37)(13,43)
\psframe[linecolor=darkgray,fillstyle=solid,fillcolor=lightgray,linewidth=0.02in](17,67)(43,93)
\psframe[linecolor=darkgray,fillstyle=solid,fillcolor=lightgray,linewidth=0.02in](47,7)(73,33)
\psframe[linecolor=darkgray,fillstyle=solid,fillcolor=lightgray,linewidth=0.02in](77,47)(93,63)
\pscircle*(10,40){0.04in}
\pscircle*(20,70){0.04in}
\pscircle*(30,90){0.04in}
\pscircle*(40,80){0.04in}
\pscircle*(50,30){0.04in}
\pscircle*(60,20){0.04in}
\pscircle*(70,10){0.04in}
\pscircle*(80,50){0.04in}
\pscircle*(90,60){0.04in}
\end{pspicture}
\end{center}
\caption{The plot of $479832156$, an inflation of $2413$.}\label{fig-479832156}
\end{figure}
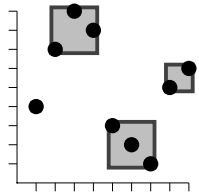

\begin{proposition}\label{simple-decomp-1}
Every permutation is the inflation of a unique simple permutation.
\end{proposition}

For simple permutations of length at least four, one can say more:

\begin{proposition}\label{simple-decomp-2}
If $\pi$ can be written as $\sigma[\alpha_1,\dots,\alpha_m]$ where $\sigma$ is simple and $m\ge 4$, then the $\alpha_i$'s are unique.
\end{proposition}

In the case where $\pi=12[\alpha_1,\alpha_2]$, some caution is needed.  A {\it sum indecomposable\/} permutation is one that cannot be written as $12[\alpha_1,\alpha_2]$, whilst a {\it skew indecomposable\/} permutation is one than cannot be written as $21[\alpha_1,\alpha_2]$.

\begin{proposition}\label{simple-decomp-3}
If $\pi$ is an inflation of $12$, then there is a unique sum indecomposable $\alpha_1$ such that $\pi=12[\alpha_1,\alpha_2]$ for some $\alpha_2$.  The same holds with $12$ replaced by $21$ and ``sum'' replaced by ``skew''.
\end{proposition}

We refer to the unique decompositions guaranteed by Propositions~\ref{simple-decomp-1}--\ref{simple-decomp-3} as the {\it substitution decomposition\/}.

Let $\Si(X)$ denote the simple permutations in the set $X$ of permutations.  The {\it wreath closure\/}, $\W(X)$, of $X$ is the largest permutation class with $\Si(\W(X))=\Si(X)$.  For example, the wreath closure of $\Av(132)$ is the largest class whose only simple permutations are $1$, $12$, and $21$.  This class is known as the {\it separable permutations\/}\footnote{The separable permutations seem to have made their first appearance as the permutations that can be sorted by pop-stacks in series, see Avis and Newborn~\cite{avis:on-pop-stacks-i:}.  Shapiro and Stephens~\cite{shapiro:bootstrap-perco:} showed that the separable permutations are those that fill up under bootstrap percolation.  The separable permutations are essentially the permutation analogue of series-parallel posets (see Stanley~\cite[Section 3.2]{stanley:enumerative-com:}) and complement reducible graphs (see Corneil, Lerchs, and Burlingham~\cite{corneil:complement-redu:}).  Their enumeration is given by the large Schr\"oder numbers (see Footnote~3 or Example~\ref{ex-separable}).}, $\Av(2413,3142)$.  When $\W(\C)=\C$, $\C$ is said to be {\it wreath-closed\/}%
\footnote{%
It is quite easy to decide if a permutation class given by a finite basis is wreath-closed:
\begin{proposition}[Atkinson and Stitt~\cite{atkinson:restricted-perm:a}]\label{decide-wreath-complete}
A permutation class is wreath-closed if and only if each of its basis elements is simple.
\end{proposition}
One may also wish to compute the basis of $\W(\C)$.  Murphy includes in his thesis~\cite{murphy:restricted-perm:} an example of a finitely based class whose wreath closure is infinitely based.  The natural question is then:
\begin{question}\label{question-wreath-closure}
Given a finite basis $B$, is it decidable whether $\W(\Av(B))$ is finitely based?
\end{question}
The analogous question for graphs was raised by Giakoumakis~\cite{giakoumakis:on-the-closure-:} and has received some attention, see for example Zverovich~\cite{zverovich:a-finiteness-th:}.  In one notable case Question~\ref{question-wreath-closure} has a nice answer:
\begin{proposition}\label{prop-wc-fin-simple}
If the longest simple permutations in $\C$ have length $k$ then the basis elements of $\W(\C)$ have length at most $k+2$.
\end{proposition}
\begin{proof}
The basis of $\W(\C)$ is easily seen to consist of the minimal (under the pattern-containment order) simple permutations not contained in $\C$.  Let $\pi$ be such a permutation of length $n$.  Theorem~\ref{thm-schmerl-trotter} shows that $\pi$ contains a simple permutation $\sigma$ of length $n-1$ or $n-2$.  If $n\ge k+3$, then $\sigma\notin\C$, so $\sigma\notin\W(\C)$ and thus $\pi$ does not lie in the basis of $\W(\C)$.
\end{proof}%
}
The crucial property of a wreath-closed class $\C$ is that $\sigma[\alpha_1,\dots,\alpha_m]\in\C$ for all $\sigma,\alpha_1,\dots,\alpha_m\in\C$.

\section{A More General Theorem}\label{sp1-properties}

\renewcommand{\P}{\mathcal{P}}
\newcommand{\Q}{\mathcal{Q}}
\newcommand{\R}{\mathcal{R}}
\renewcommand{\S}{\mathcal{S}}

A {\it property\/}, $P$, is any set of permutations.  We say that $\pi$ {\it satisfies\/} $P$ if $\pi\in P$.  Given a permutation class $\C$ and set $\P$ of properties, we write $\C_\P$ for the set of permutations in $\C$ that satisfy every property in $\P$, and write $f_\P$ for the generating function of $\C_\P$.

A set $\P$ of properties is said to be {\it query-complete\/} if, for all permutations $\sigma\in S_m$ and $P\in\P$, it can be decided whether $\sigma[\alpha_1,\dots,\alpha_m]$ satisfies $P$ given only the knowledge about whether $\alpha_i$ satisfies $Q$, for all $i\in[m]$ and $Q\in\P$.  In other words, this condition means that every query ``does $\sigma[\alpha_1,\dots,\alpha_m]$ satisfy $P$?'' can be answered by considering only questions of the form ``does $\alpha_i$ satisfy $Q$?''  Note that we need only check the condition for nontrivial simple $\sigma$'s, as $\sigma[\alpha_1,\dots,\alpha_m]$ can always be expressed as the inflation of a simple permutation (by Proposition~\ref{simple-decomp-1}).

For example, $\{\Av(132)\}$ is not query-complete, as witnessed by the fact that $12[1,1]\in\Av(132)$ but $12[1,21]\notin\Av(132)$.  However, $\{\Av(132),\Av(21)\}$ is query-complete:
\begin{eqnarray*}
12[\alpha_1,\alpha_2]\in\Av(132)&\iff&\alpha_1\in\Av(132)\mbox{ and }\alpha_2\in\Av(21),\\
21[\alpha_1,\alpha_2]\in\Av(132)&\iff&\alpha_1\in\Av(132)\mbox{ and }\alpha_2\in\Av(132),\\
\sigma[\alpha_1,\dots,\alpha_m]\notin\Av(132)&\mbox{if}&\sigma\notin\{1,12,21\},\\
12[\alpha_1,\alpha_2]\in\Av(21)&\iff&\alpha_1\in\Av(21)\mbox{ and }\alpha_2\in\Av(21),\\
\sigma[\alpha_1,\dots,\alpha_m]\notin\Av(21)&\mbox{if}&\sigma\notin\{1,12\},
\end{eqnarray*}
and since $12$ and $21$ are the only nontrivial simple permutations in $\Av(132)$, these are the only computations that need to be performed.

Our main result in this section is the following.

\begin{theorem}\label{sp1-really-main}
Let $\C$ be a wreath-closed permutation class containing only finitely many simple permutations, $\P$ a finite query-complete set of properties, and $\Q\subseteq\P$.  The generating function for the set of permutations in $\C$ satisfying every property in $\Q$, i.e., $f_\Q$, is algebraic.
\end{theorem}

Let us postpone the proof of Theorem~\ref{sp1-really-main} momentarily in order to introduce notation and review pertinent facts.  We begin by considering the case $\P=\emptyset$, which contains many of the main ideas of the proof in a more digestible form.  We introduce two properties,
\begin{eqnarray*}
\noplus&=&\{\mbox{sum indecomposable permutations}\}\mbox{ and}\\
\nominus&=&\{\mbox{skew indecomposable permutations}\}.
\end{eqnarray*}
Note that both $\{\noplus\}$ and $\{\nominus\}$ are query-complete, because for simple $\sigma$,
$$
\sigma[\alpha_1,\dots,\alpha_m]\in\noplus\iff\sigma\neq 12.
$$
We also introduce the notation
$$
\sigma[\C_1,\dots,\C_m]=\{\sigma[\alpha_1,\dots,\alpha_m] : \mbox{$\alpha_i\in\C_i$ for all $i\in[m]$}\}.
$$

By Propositions~\ref{simple-decomp-1}--\ref{simple-decomp-3} and the assumption that $\C$ is wreath-closed, $\C$ can be written as
$$
\C=\{1\}\uplus 12[\C_{\noplus},\C]\uplus 21[\C_{\nominus},\C]\uplus\biguplus_{\substack{\sigma\in\Si(\C)\\|\sigma|\ge 4}} \sigma[\C,\dots,\C],
$$
while $\C_{\noplus}$ and $\C_{\nominus}$ have the expressions
$$
\begin{array}{rclcl}
\C_{\noplus}&=&\{1\}\uplus 21[\C_{\nominus},\C]\uplus
\displaystyle\biguplus_{\substack{\sigma\in\Si(\C)\\|\sigma|\ge 4}}\sigma[\C,\dots,\C]
&=&
\C\setminus 12[\C_\noplus,\C],
\\
\C_{\nominus}&=&\{1\}\uplus 12[\C_{\noplus},\C]\uplus
\displaystyle\biguplus_{\substack{\sigma\in\Si(\C)\\|\sigma|\ge 4}} \sigma[\C,\dots,\C]
&=&
\C\setminus 21[\C_\nominus,\C].
\end{array}
$$
These give the system
$$
\left\{
\begin{array}{rclcl}
f&=&x+f_{\noplus}f+f_{\nominus}f+\displaystyle\sum_{\substack{\sigma\in\Si(\C)\\|\sigma|\ge 4}} f^{|\sigma|},\\
f_{\noplus}&=&x+f_{\nominus}f+\displaystyle\sum_{\substack{\sigma\in\Si(\C)\\|\sigma|\ge 4}} f^{|\sigma|}
&=&\displaystyle\frac{f}{1+f},\\
f_{\nominus}&=&x+f_{\noplus}f+\displaystyle\sum_{\substack{\sigma\in\Si(\C)\\|\sigma|\ge 4}} f^{|\sigma|}
&=&\displaystyle\frac{f}{1+f}.
\end{array}
\right.
$$
If we now let $s$ denote the generating function for the simple permutations of length at least $4$ in $\C$, we get that
$$
f=x+\frac{2f^2}{1+f}+s(f),
$$
so if $s$ is algebraic, a fortiori if $s$ is polynomial, $f$ is algebraic%
\footnote{In particular, note that the separable permutations correspond to $s=0$; making this substitution leaves $f=x+2f^2/(1+f)$, giving the large Schr\"oder numbers.%
\newcounter{fn-separable-comment}\setcounter{fn-separable-comment}{\value{footnote}}}.

In general the situation will not be so straightforward and we will resort to algebraic systems.  The following is a specialisation of the more general treatment in Stanley~\cite[Section 6.6]{stanley:enumerative-com:a}.

Let $A=\{a_1,\dots,a_n\}$ denote an alphabet.  A {\it proper algebraic system\/} over $\mathbb{Q}[x_1,\dots,x_m]$ is a set of equations $a_i=p_i(x_1,\dots,x_m,a_1,\dots,a_n)$ where each $p_i$ is a polynomial with coefficients from $\mathbb{Q}$, has constant term $0$, and contains no terms of the form $ca_i$ where $c\in\mathbb{Q}$.  The solution to such a system is a tuple $(f_1,\dots,f_n)$ of formal power series from $\mathbb{Q}[[x_1,\dots,x_m]]$ such that  for all $i$, $f_i$ is equal to $p_i(x_1,\dots,x_m,a_1,\dots,a_n)$ evaluated at $(a_1,\dots,a_n)=(f_1,\dots,f_n)$.

\begin{theorem}[Stanley~{\cite[Proposition 6.6.3 and Theorem 6.6.10]{stanley:enumerative-com:a}}]\label{prop-alg-system}
Every proper algebraic system $(p_1,\dots,p_n)$ over $\mathbb{Q}[x_1,\dots,x_m]$ has a unique solution $(f_1,\dots,f_n)$.  Moreover, each of these $f_i$'s is algebraic over $\mathbb{Q}[x_1,\dots,x_m]$.
\end{theorem}

\bigskip\noindent{\it Proof of Theorem~\ref{sp1-really-main}.\/}  If $\P$ is query-complete then $\P\cup\{\noplus,\nominus\}$ is also query-complete, so we may assume without loss that $\noplus,\nominus\in\P$.  

Let $\P(\pi)$ denote the set of properties in $\P$ satisfied by $\pi$ and let $g_\R$ denote the generating function for the set of $\pi\in\C$ with $\P(\pi)=\R$, so
$$
f_\Q=\sum_{\Q\subseteq\R\subseteq\P} g_\R.
$$
As $\P$ is query-complete, for each simple $\sigma\in\C$ of length $m$ there is a finite collection of $m$-tuples of sets of properties such that $\P(\sigma[\alpha_1,\dots,\alpha_m])=\R$ precisely if $(\P(\alpha_1),\dots,\P(\alpha_m))$ lies in this collection.  If $m\ge 4$ then Propositions~\ref{simple-decomp-1} and \ref{simple-decomp-2} imply that the generating function for all inflations $\pi$ of $\sigma$ with $\P(\pi)=\R$ can be expressed nontrivially as a polynomial in $\{g_\S : \S\subseteq\P\}$.  If $m=2$, suppose $\sigma=12$ without loss.  By Proposition~\ref{simple-decomp-3}, all inflations of $12$ have a unique decomposition as $12[\alpha_1,\alpha_2]$ where $\alpha_1\in\noplus$.  Thus the generating function for inflations $\pi$ of $12$ with $\P(\pi)=\R$ can be expressed as a sum of terms of the form $g_\S g_\mathcal{T}$ where $\noplus\in\S$.

Therefore $g_\R$ can be expressed as a polynomial in $x$ (depending on whether $\P(1)=\R$) and $\{g_\S:\S\subseteq\P\}$.  Moreover, these polynomials have no constant terms and no terms of the form $cg_\S$ for constant $c\neq 0$.  Thus they form a proper algebraic system, so Theorem~\ref{prop-alg-system} implies that each $g_\S$ is algebraic.
\qed\bigskip

We are now in a position to indicate the proof of Theorem~\ref{sp1-main}, which will be completed over the course of the next two sections.  Suppose that $\C$ contains only finitely many simple permutations and let $B$ denote the basis of $\C$.  A result of Higman~\cite{higman:ordering-by-div:} implies that $\C$ is finitely based, i.e., $B$ is finite; see Albert and Atkinson~\cite{albert:simple-permutat:} for details.  By the forthcoming Lemma~\ref{query-complete-beta}, each property $\Av(\beta)$ lies in a finite query-complete set.  Clearly the union over all $\beta\in B$ of these finite query-complete sets is again a finite query-complete set, and we have
$$
\C=\W(\C)_{\{\Av(\beta) : \beta\in B\}},
$$
so the generating function for $\C$ is algebraic by Theorem~\ref{sp1-really-main}.  Now suppose that we are interested in the even permutations in $\C$.  In this case we let $EV$ denote the set of even permutations.  Lemma~\ref{query-complete-even} shows that $EV$ lies in a finite query-complete set, so the even permutations in $\C$,
$$
\C_{\{EV\}}=\W(\C)_{\{EV\}\cup\{\Av(\beta) : \beta\in B\}},
$$
also have an algebraic generating function by Theorem~\ref{sp1-really-main}, thereby establishing that part of Theorem~\ref{sp1-main}.  The other cases --- with the exception of involutions, discussed in Section~\ref{sp1-involution} --- are analogous.

\section{Finite Query-Complete Sets}\label{sp1-fqcs-exs}

In order to utilise the general framework provided by Theorem~\ref{sp1-really-main}, we need to construct finite query-complete sets of properties.  We do that here, beginning with the most important one.

\begin{lemma}\label{query-complete-beta}
For every permutation $\beta$, the set $\{\Av(\delta) : \delta\le\beta\}$ is query-complete.
\end{lemma}
\begin{proof}
By applying the lemma to $\{\Av(\gamma) : \gamma\le\delta\}$, we see that it suffices to prove that whether $\pi=\sigma[\alpha_1,\dots,\alpha_m]\in\Av(\beta)$ can be decided entirely by knowing, for each $i$, which permutations $\delta$ satisfy $\delta\le\alpha_i$ and $\delta\le\beta$.

We define a {\it lenient inflation\/} to be an inflation $\sigma[\alpha_1,\dots,\alpha_m]$ in which the $\alpha_i$'s are allowed to be empty.  List all expressions of $\beta$ as a lenient inflation of $\sigma$ as
\begin{eqnarray*}
\beta&=&\sigma[\beta_1^{(1)},\dots,\beta_m^{(1)}],\\
&\vdots&\\
\beta&=&\sigma[\beta_1^{(t)},\dots,\beta_m^{(t)}].
\end{eqnarray*}
Clearly if we have, for some $s\in[t]$, $\alpha_i\ge\beta_i^{(s)}$ for all $i\in[m]$, then $\pi\ge\beta$.  Equivalently, to have $\pi\in\Av(\beta)$, for every $s\in[t]$ there must be at least one $i\in[m]$ for which $\alpha_i\not\ge\beta_i^{(s)}$.  Conversely, every embedding of $\beta$ into $\pi$ gives one of the lenient inflations in the list above, completing the proof.
\end{proof}

The proof of Lemma~\ref{query-complete-beta} extends in the obvious manner to show that the property of avoiding a blocked or barred permutation (or, for that matter, a permutation combining these restrictions) also lies in a finite query-complete set.

\begin{lemma}\label{query-complete-alternating}
The set of properties consisting of
\begin{itemize}
\item $AL=\{\mbox{alternating permutations}\}$,
\item $BR=\{\mbox{permutations beginning with a rise, i.e., permutations with $\pi(1)<\pi(2)$}\}$,
\item $ER=\{\mbox{permutations ending with a rise}\}$, and
\item $\{1\}$.
\end{itemize}
is query-complete.
\end{lemma}
\begin{proof}
Clearly $\{\{1\},BR,ER\}$ is query-complete:
\begin{eqnarray*}
\sigma[\alpha_1,\dots,\alpha_m]\in BR&\iff&
\alpha_1\in BR
\mbox{ or }
\left(\alpha_1=1\mbox{ and }\sigma(1)<\sigma(2)\right),
\\
\sigma[\alpha_1,\dots,\alpha_m]\in ER&\iff&
\alpha_m\in ER
\mbox{ or }
\left(\alpha_m=1\mbox{ and }\sigma(m-1)<\sigma(m)\right).
\end{eqnarray*}

For $\sigma[\alpha_1,\dots,\alpha_m]$ to be an alternating permutation, we first need $\alpha_1,\dots,\alpha_m\in AL$.  Now note that if $\sigma(i)<\sigma(i+1)$ then $\sigma[\alpha_1,\dots,\alpha_m]$ will contain a rise between the interval corresponding to $\sigma(i)$ and the interval corresponding to $\sigma(i+1)$.  Thus $\alpha_i$ must end --- and $\alpha_{i+1}$ must begin --- with a descent (or have length $1$), so we need $\alpha_i\notin ER$ and $\alpha_{i+1}\notin BR$.  Similarly, if $\sigma(i)>\sigma(i+1)$ then we need $\alpha_i\in ER\cup\{1\}$ and $\alpha_{i+1}\in BR\cup\{1\}$.  These are easily seen to also be sufficient conditions for $\sigma[\alpha_1,\dots,\alpha_m]$ to be an alternating permutation, completing the proof.
\end{proof}

\begin{lemma}\label{query-complete-even}
The set of properties consisting of
\begin{itemize}
\item $EV=\{\mbox{even permutations}\}$ and
\item $EL=\{\mbox{permutations of even length}\}$
\end{itemize}
is query-complete.
\end{lemma}
\begin{proof}
We have
$$
\sigma[\alpha_1,\dots,\alpha_m]\in EL
\iff
\mbox{an even number of $\alpha_i$'s fail to lie in $EL$},
$$
so $\{EL\}$ is query-complete.  To see that $\{EV,EL\}$ is query-complete, we divide the inversions in $\sigma[\alpha_1,\dots,\alpha_m]$ into two groups: inversions within a single block $\alpha_i$ and inversions between two blocks $\alpha_i$ and $\alpha_j$.  We need to compute the parity of each of these numbers.  The parity of the number of inversions within $\alpha_i$ is given by whether $\alpha_i\in EV$.  The parity of the number of inversions between the blocks $\alpha_i$ and $\alpha_j$, $i<j$, is even if $\sigma(i)<\sigma(j)$, if $\alpha_i\in EL$, or if $\alpha_j\in EL$, and odd otherwise.
\end{proof}

\begin{lemma}\label{query-complete-Dumont}
The set of properties consisting of
\begin{itemize}
\item $DU=\{\mbox{Dumont permutations of the first kind}\}$ and
\item $EL=\{\mbox{permutations of even length}\}$
\end{itemize}
is query-complete.
\end{lemma}
\begin{proof}
It suffices to determine which entries of $\sigma[\alpha_1,\dots,\alpha_m]$ have even value and which have odd value, and this can be decided based on the knowledge of which $\alpha_i$'s have even length.
\end{proof}

\section{Involutions}\label{sp1-involution}

We have now proved all of Theorem~\ref{sp1-main} except for the claim that we can count the involutions in $\C$ (and even involutions, alternating involutions, etc.).  Unfortunately, involutionhood lies just outside the scope of our query-complete-property machinery: letting $I$ denote the set of involutions we have that $12[\alpha_1,\alpha_2]\in I\iff \alpha_1,\alpha_2\in I$, but when is $21[\alpha_1,\alpha_2]\in I$?

We begin by considering the effect of inversion on the substitution decomposition.  First note that if $\pi$ is an inflation of $\sigma$ then $\pi^{-1}$ is an inflation of $\sigma^{-1}$.  Recalling Proposition~\ref{simple-decomp-1} (``every permutation is the inflation of a unique simple permutation''), we have that if $\pi$ is an involution then it must be the inflation of a simple involution.  By Proposition~\ref{simple-decomp-2} we then get the following:

\begin{proposition}\label{involution-decomp-1}
If $\pi=\sigma[\alpha_1,\dots,\alpha_m]$ is an involution and $\sigma$ is a simple permutation of length at least $4$, then $\sigma$ is an involution and $\alpha_i=\alpha_{\sigma(i)}^{-1}$ for all $i\in[m]$.
\end{proposition}

As in Section~\ref{sp1-inflations}, $\sigma=12$ and $\sigma=21$ must be handled separately.  Moreover, when $\pi=21[\alpha_1,\alpha_2]$, unless $\alpha_1$ and $\alpha_2$ are both skew indecomposable, we use the {\it middle greedy decomposition\/}, expressing $\pi$ as $321[\alpha_1,\alpha_2,\alpha_3]$ with $\alpha_2$ as long as possible.

\begin{proposition}\label{involution-decomp-2}
The involutions that are inflations of permutations of length $2$ are precisely those of the form
\begin{itemize}
\item $12[\alpha_1,\alpha_2]$ for involutions $\alpha_1$ and $\alpha_2$,
\item $21[\alpha_1,\alpha_2]$ for skew indecomposable $\alpha_1$ and $\alpha_2$ with $\alpha_1=\alpha_2^{-1}$, and
\item $321[\alpha_1,\alpha_2,\alpha_3]$, where $\alpha_1$ and $\alpha_3$ are skew indecomposable, $\alpha_1=\alpha_3^{-1}$, and $\alpha_2$ is an involution (i.e., the middle greedy decompositions).
\end{itemize}
\end{proposition}

Define the {\it inverse\/} of the property $P$ by $P^{-1}=\{\pi^{-1} : \pi\in P\}$, and for a set of properties $\P$, $\P^{-1}=\{P^{-1} : P\in\P\}$.
%
%

\begin{theorem}\label{sp1-involutions}
Let $\C$ be a wreath-closed permutation class containing only finitely many simple permutations, $\P$ a finite query-complete set of properties, and $\Q\subseteq\P$.  The generating function for the set of {\it involutions\/} in $\C$ satisfying every property in $\Q$ is algebraic.
\end{theorem}
\begin{proof}
We assume (without loss) both that $\noplus,\nominus\in\P$ and that $\P=\P^{-1}$.  As in the proof of Theorem~\ref{sp1-really-main}, let $\P(\pi)$ denote the set of properties in $\P$ satisfied by $\pi$ and $g_\R$ denote the generating function for the set of $\pi\in\C$ with $\P(\pi)=\R$.  Also let $h_\R$ denote the generating function for the set of involutions $\pi\in\C$ with $\P(\pi)=\R$.  It suffices to show that each $h_\R$ is algebraic over $\mathbb{Q}[x]$.

As Propositions~\ref{involution-decomp-1} and \ref{involution-decomp-2} indicate, we need to count pairs $(\alpha,\alpha^{-1})$ where $\alpha$ and $\alpha^{-1}$ satisfy certain sets of properties.  To this end define
$$
p_\R = \sum_{\substack{\alpha\in\C\\\P(\alpha)=\R}} x^{|\alpha|+|\alpha^{-1}|}.
$$
(Note that if $\P(\alpha)=\R$ then $\P(\alpha^{-1})=\R^{-1}$ because $\P=\P^{-1}$.)

Now take $\sigma$ to be a simple permutation.  We need to compute the contribution to $h_\R$ of inflations of $\sigma$.  As before, since $\P$ is query-complete, $\P(\sigma[\alpha_1,\dots,\alpha_m])=\R$ if and only if $(\P(\alpha_1),\dots,\P(\alpha_m))$ lies in a certain collection of $m$-tuples of sets of properties.  Choose one of these $m$-tuples, say $(\R_1,\dots,\R_m)$, and suppose first that $m=|\sigma|\ge 4$.  It suffices to calculate the contribution of involutions of the form $\sigma[\alpha_1,\dots,\alpha_m]$ with $\P(\alpha_i)=\R_i$ for all $i\in[m]$.  If there is some $j\in[m]$ for which $\R_j\neq \R_{\sigma(j)}^{-1}$ then this contribution is $0$.  Otherwise the contribution is a single term in which each fixed point $j$ corresponds to an $h_{\R_j}$ factor and each non-fixed-point pair $(j,\sigma(j))$ corresponds to a $p_{\R_j}$ factor.  A similar analysis of inflations of $12$ and $21$ --- in the latter case using the middle greedy decomposition when necessary --- allows us to compute the contribution of these inflations.

Therefore each $h_\R$ can be expressed nontrivially as a polynomial in $x$, $\{h_\S : \S\subseteq\P\}$, and $\{p_\S : \S\subseteq\P\}$.  Viewing $x$ and $\{p_\S : \S\subseteq\P\}$ as variables, Theorem~\ref{prop-alg-system} implies that each $h_\R$ is algebraic over $\mathbb{Q}[x,\{p_\S:\S\subseteq\P\}]$.  However, $p_\S$ is nothing other than $g_\S(x^2)$, so $\mathbb{Q}[x,\{p_\S:\S\subseteq\P\}]$ is an algebraic extension of $\mathbb{Q}[x]$ by Theorem~\ref{sp1-really-main}, proving the theorem.
\end{proof}

One could adapt Theorem~\ref{sp1-involutions} to count the permutations in $\C$ that are invariant under other symmetries.  For example, permutations invariant under the composition of reverse and complement have received some attention, see Guibert and Pergola~\cite{gp:vex}.

\section{Examples}\label{sp1-examples}

While we have already shown how to enumerate separable permutations in Footnote~\arabic{fn-separable-comment}, here we use the approach of Theorem~\ref{sp1-really-main}.

\begin{example}[Separable permutations]\label{ex-separable}
With the notation from the proof of Theorem~\ref{sp1-really-main}, we have that for the separable permutations:
$$
\left\{\begin{array}{lcl}
g_{\noplus,\nominus}&=&x,\\
g_\noplus&=&(g_{\noplus,\nominus}+g_\nominus)(g_{\noplus,\nominus}+g_\noplus+g_\nominus),\\
g_\nominus&=&(g_{\noplus,\nominus}+g_\noplus)(g_{\noplus,\nominus}+g_\noplus+g_\nominus).
\end{array}\right.
$$
We are interested in $f=g_{\noplus,\nominus}+g_\noplus+g_\nominus$.  We can find the algebraic equation satisfied by $f$ by computing the Gr\"obner basis (under an appropriate monomial ordering) of the ideal generated by
$$
\begin{array}{l}
g_{\noplus,\nominus}-x,\\
g_\noplus-(g_{\noplus,\nominus}+g_\nominus)(g_{\noplus,\nominus}+g_\noplus+g_\nominus),\\
g_\nominus-(g_{\noplus,\nominus}+g_\noplus)(g_{\noplus,\nominus}+g_\noplus+g_\nominus),\\
f-g_{\noplus,\nominus}+g_\noplus+g_\nominus.
\end{array}
$$
This calculation shows that
$$
f^2+(x-1)f+x=0,
$$
so
$$
f=\frac{1-x-\sqrt{1-6x+x^2}}{2}.
$$
This, reassuringly, is the generating function for the large Schr\"oder numbers.
\end{example}

This system does not change dramatically when another simple permutation is introduced, as shown by the next example.

\begin{example}[The wreath closure of $\mathbf{1, 12, 21,}$ and $\mathbf{2413}$%
\footnote{Using Proposition~\ref{prop-wc-fin-simple} it can be computed that this is the class $\Av(3142,25314,246135,362514)$.}%
]\label{ex-wc-12-21-2413}
Here the system is
$$
\left\{\begin{array}{lcl}
g_{\noplus,\nominus}&=&x+(g_{\noplus,\nominus}+g_\noplus+g_\nominus)^4,\\
g_\noplus&=&(g_{\noplus,\nominus}+g_\nominus)(g_{\noplus,\nominus}+g_\noplus+g_\nominus),\\
g_\nominus&=&(g_{\noplus,\nominus}+g_\noplus)(g_{\noplus,\nominus}+g_\noplus+g_\nominus).\\
\end{array}\right.
$$
Setting
$$
f=g_{\noplus,\nominus}+g_\noplus+g_\nominus
$$
and computing the Gr\"obner basis for the associated ideal as in the previous example shows that the generating function for this class satisfies
$$
f^5+f^4+f^2+(x-1)f+x=0.
$$
\end{example}

\begin{example}[Av($\mathbf{132}$)]\label{ex-count-132}
The wreath closure of $\Av(132)$ is the class of separable permutations, so to enumerate $\Av(132)$ we need to refine Example~\ref{ex-separable}.  While Proposition~\ref{query-complete-beta} shows that $\{\Av(1),\Av(12),\Av(21),\Av(132)\}$ is query-complete, we remarked at the beginning of Section~\ref{sp1-properties} that $\{\Av(21),\Av(132)\}$ will suffice.  Our system is then
$$
\left\{\begin{array}{lcl}
g_{\noplus,\nominus,\Av(21)}&=&x,\\
g_{\nominus,\Av(21)}
&=&
g_{\nominus,\noplus,\Av(21)}(g_{\noplus,\nominus,\Av(21)}+g_{\nominus,\Av(21)}),\\
g_{\noplus}
&=&
(g_{\noplus,\nominus,\Av(21)}+g_{\nominus,\Av(21)}+g_{\nominus})
(g_{\noplus,\nominus,\Av(21)}+g_{\nominus,\Av(21)}+g_{\noplus}+g_{\nominus}),\\
g_{\nominus}
&=&
g_{\noplus}(g_{\noplus,\nominus,\Av(21)}+g_{\nominus,\Av(21)}).
\end{array}\right.
$$
(As we are only interested in $132$-avoiding permutations we have suppressed the subscript $\Av(132)$, which would otherwise be present in all these terms.)  Setting
$$
f=g_{\noplus,\nominus,\Av(21)}+g_{\nominus,\Av(21)}+g_\noplus+g_\nominus
$$
and computing a Gr\"obner basis as before yields
$$
f=\frac{1-2x-\sqrt{1-4x}}{2x},
$$
the generating function for the Catalan numbers, as expected.
\end{example}

\begin{example}[Av($\mathbf{2413, 3142, 2143}$)]\label{ex-sep-2143}
Our system in this case is
$$
\left\{\begin{array}{lcl}
g_{\noplus,\nominus,\Av(21)}
&=&
x,
\\
g_{\nominus,\Av(21)}
&=&
g_{\noplus,\nominus,\Av(21)}
(g_{\noplus,\nominus,\Av(21)}+g_{\nominus,\Av(21)}),
\\
g_{\noplus}
&=&
(g_{\noplus,\nominus,\Av(21)}+g_{\nominus,\Av(21)}+g_{\nominus})
(g_{\noplus,\nominus,\Av(21)}+g_{\nominus,\Av(21)}+g_{\noplus}+g_{\nominus}),
\\
g_{\nominus}
&=&
g_{\noplus,\nominus,\Av(21)}
(g_{\noplus}+g_{\nominus})
+
g_{\noplus}(g_{\noplus,\nominus,\Av(21)}+g_{\nominus,\Av(21)}),
\end{array}\right.
$$
where here we have suppressed the $\Av(2143)$ subscript.  Solving as before shows that the generating function for this class is
$$
\frac{1-3x+2x^2-\sqrt{1-6x+5x^2}}{2x(2-x)},
$$
and thus the number of permutations of length $n$ in this class is $\sum{n\choose k}f_{n-k}$ where $(f_n)$ denotes Fine's sequence.
\end{example}

\begin{example}[Alternating separable permutations]\label{ex-alt-separable}
Lemma~\ref{query-complete-alternating} shows that we need to introduce the properties $AL$ (alternating permutations), $BR$ (permutations beginning with a rise), $ER$ (permutations ending with a rise), and $\{1\}$.  In the separable case $\{1\}=\noplus\cap\nominus$ so we can ignore that property, and as $AL$ occurs in each of the terms of our system we suppress it.  We then have
$$
\left\{
\begin{array}{lcl}
g_{\noplus,\nominus}&=&x,\\
g_{\noplus}&=&(g_{\noplus,\nominus}+g_{\nominus,ER})(g_{\noplus,\nominus}+g_{\noplus,BR}+g_{\nominus,BR}),\\
g_{\noplus,BR}&=&g_{\nominus,BR,ER}(g_{\noplus,\nominus}+g_{\noplus,BR}+g_{\nominus,BR}),\\
g_{\noplus,ER}&=&(g_{\noplus,\nominus}+g_{\nominus,ER})(g_{\noplus,BR,ER}+g_{\nominus,BR,ER}),\\
g_{\noplus,BR,ER}&=&g_{\nominus,BR,ER}(g_{\noplus,BR,ER}+g_{\nominus,BR,ER}),\\
g_{\nominus}&=&g_{\noplus}(g_{\noplus}+g_{\nominus}),\\
g_{\nominus,BR}&=&(g_{\noplus,\nominus}+g_{\noplus,BR})(g_{\noplus}+g_{\nominus}),\\
g_{\nominus,ER}&=&g_{\noplus}(g_{\noplus,\nominus}+g_{\noplus,ER}+g_{\nominus,ER}),\\
g_{\nominus,BR,ER}&=&(g_{\noplus,\nominus}+g_{\noplus,BR})(g_{\noplus,\nominus}+g_{\noplus,ER}+g_{\nominus,ER}).
\end{array}\right.
$$
The generating function for these permutations satisfies
$$
f^3-(2x^2-5x+4)f^2-(4x^3+x^2-8x)f-(2x^4+5x^3+4x^2)=0.
$$
\end{example}

\begin{example}[Separable involutions]\label{ex-separable-involution}
Using the notation from the proof of Theorem~\ref{sp1-involutions}, we wish to find $f=h_{\noplus, \nominus}+h_\noplus+h_\nominus$.  These generating functions are related to each other and to the $p$ generating functions by
$$
\left\{
\begin{array}{lcl}
h_{\noplus,\nominus}&=&x,\\
h_\noplus&=&(p_{\noplus,\nominus}+p_{\nominus})+(p_{\noplus,\nominus}+p_{\nominus})(h_{\noplus,\nominus}+h_\noplus+h_\nominus),\\
h_\nominus&=&(h_{\noplus,\nominus}+h_\noplus)(h_{\noplus,\nominus}+h_\noplus+h_\nominus).
\end{array}
\right.
$$
From Example~\ref{ex-separable} it can be computed that
\begin{eqnarray*}
p_{\noplus,\nominus}-x^2&=&0,\\
2p_\noplus^2+(3x^2-1)p_\noplus+x^4&=&0,\\
2p_\nominus^2+(3x^2-1)p_\nominus +x^4&=&0.
\end{eqnarray*}
Combining these with the system above and solving as usual shows that
$$
x^2f^4 + (x^3+3x^2+x-1)f^3 + (3x^3+6x^2-x)f^2  + (3x^3+7x^2-x-1)f +x^3+3x^2+x=0.
$$
\end{example}


%
%
%
%
%
%
%
%
%
%
%
%
%
%
%
%

\section{Concluding Remarks}

We end by discussing issues related to the applicability and application of these techniques.

\minisec{Determining if these methods apply}
The results we have established apply only to permutation classes with finitely many simple permutations.  Thus it would be useful to be able to determine whether a permutation class contains finitely many simple permutations.  This can be done:

\begin{theorem}[Brignall, Ru\v{s}kuc, and Vatter~\cite{brignall:simple-permutat:b}]
It is decidable whether permutation class given by a finite basis contains only finitely many simple permutations.
\end{theorem}

\minisec{Finding the simple permutations}
Thus far we have tacitly assumed that the set of simple permutations in our class is known.  Since classes are often specified by their bases, this set of simple permutations must first be computed.  Assuming that this set is finite, it can be computed via a result of Schmerl and Trotter.  While we state only the permutation case (a proof of this case is also given by Murphy~\cite{murphy:restricted-perm:}), their result covers all irreflexive binary relational structures.  See Ehrenfeucht and McConnell~\cite{ehrenfeucht:a-k-structure-g:} for a version of this theorem for certain other relational structures.

\begin{theorem}[Schmerl and Trotter~\cite{schmerl:critically-inde:}]\label{thm-schmerl-trotter}
Every simple permutation of length $n\ge 2$ contains a simple permutation of length $n-1$ or $n-2$.
\end{theorem}

For example, the number of simple permutations in $\Av(1324,2143,4231)$ of lengths $1$ to $7$ is $1,2,0,2,4,0,0$.  Because there are no simple permutations of length $6$ or $7$ in this class, Theorem~\ref{thm-schmerl-trotter} ensures that it contains no longer simple permutations.

\minisec{Derangements}
Notably absent from Theorem~\ref{sp1-main} is the number of derangements in $\C$, despite the fact that $132$-avoiding derangements are counted by Fine's sequence (Robertson, Saracino, and Zeilberger~\cite{robertson:refined-restric:}).  To see that the set of derangements does not lie in a finite query-complete set of properties, for $\alpha\in S_n$ define $D(\alpha)=\{\alpha(i)-i:i\in[n]\}$.  Then $21[12\cdots j,\alpha]$ is a derangement if and only if $j\notin D(\alpha)$.  This shows that $\alpha_1$ and $\alpha_2$ must lie in different sets of properties whenever $D(\alpha_1)\cap\mathbb{N}\neq D(\alpha_2)\cap\mathbb{N}$, implying that the set of derangements can only lie in an infinite query-complete set of properties.  

\minisec{Other structures}
It is natural to expect the substitution decomposition to lead to enumerative results for other types of combinatorial object.  Indeed, Stanley~\cite{stanley:enumeration-of-:} used a method similar to the substitution decomposition to count unlabelled {\it series-parallel posets\/}, those generated from the singleton poset $1$ by ordinal sum and disjoint union.  While in the permutation case we are fortunate to be able to ignore isomorphism concerns, this is not the case for series-parallel posets, which do not have an algebraic generating function (Stanley gives a functional equation satisfied by their generating function).

\minisec{Acknowledgements}
We thank Mike Atkinson, Mireille Bousquet-M\'elou, and Steve Linton for their helpful comments and suggestions, and John McDermott for his technical expertise.

\bibliographystyle{acm}
\bibliography{../refs}

\def\cprime{$'$}
\begin{thebibliography}{10}

\bibitem{albert:simple-permutat:}
{\sc Albert, M.~H., and Atkinson, M.~D.}
\newblock Simple permutations and pattern restricted permutations.
\newblock {\em Discrete Math. 300}, 1-3 (2005), 1--15.

\bibitem{atkinson:restricted-perm:a}
{\sc Atkinson, M.~D., and Stitt, T.}
\newblock Restricted permutations and the wreath product.
\newblock {\em Discrete Math. 259}, 1-3 (2002), 19--36.

\bibitem{avis:on-pop-stacks-i:}
{\sc Avis, D., and Newborn, M.}
\newblock On pop-stacks in series.
\newblock {\em Utilitas Math. 19\/} (1981), 129--140.

\bibitem{babson:generalized-per:}
{\sc Babson, E., and Steingr{\'{\i}}msson, E.}
\newblock Generalized permutation patterns and a classification of the
  {M}ahonian statistics.
\newblock {\em S\'em. Lothar. Combin. 44\/} (2000), Art. B44b, 18 pp.
  (electronic).

\bibitem{bona:the-number-of-p:}
{\sc B{\'o}na, M.}
\newblock The number of permutations with exactly {$r$} {$132$}-subsequences is
  {$P$}-recursive in the size!
\newblock {\em Adv. in Appl. Math. 18}, 4 (1997), 510--522.

\bibitem{bousquet-melou:forest-link-per:}
{\sc Bousquet-M{\'e}lou, M., and Butler, S.}
\newblock Forest-like permutations.
\newblock \texttt{arxiv:math.CO/0603617}.

\bibitem{brignall:simple-permutat:a}
{\sc Brignall, R., Huczynska, S., and Vatter, V.}
\newblock Decomposing simple permutations, with enumerative consequences.
\newblock \texttt{arXiv:math.CO/0606186}.

\bibitem{brignall:simple-permutat:b}
{\sc Brignall, R., Ru\v{s}kuc, N., and Vatter, V.}
\newblock Simple permutations: decidability and unavoidable substructures.
\newblock In preparation.

\bibitem{corneil:complement-redu:}
{\sc Corneil, D.~G., Lerchs, H., and Burlingham, L.~S.}
\newblock Complement reducible graphs.
\newblock {\em Discrete Appl. Math. 3}, 3 (1981), 163--174.

\bibitem{dumont:interpretations:}
{\sc Dumont, D.}
\newblock Interpr{\'e}tations combinatoires des nombres de {G}enocchi.
\newblock {\em Duke Math. J. 41\/} (1974), 305--318.

\bibitem{egge:132-avoiding-tw:}
{\sc Egge, E.~S., and Mansour, T.}
\newblock $132$-avoiding two-stack sortable permutations, fibonacci numbers,
  and pell numbers.
\newblock {\em Discrete Appl. Math. 143}, 1-3 (2004), 72--83.

\bibitem{egge:231-avoiding-in:}
{\sc Egge, E.~S., and Mansour, T.}
\newblock $231$-avoiding involutions and {F}ibonacci numbers.
\newblock {\em Australas. J. Combin. 30\/} (2004), 75--84.

\bibitem{ehrenfeucht:a-k-structure-g:}
{\sc Ehrenfeucht, A., and McConnell, R.}
\newblock A {$k$}-structure generalization of the theory of {$2$}-structures.
\newblock {\em Theoret. Comput. Sci. 132}, 1-2 (1994), 209--227.

\bibitem{elizalde:restricted-motz:}
{\sc Elizalde, S., and Mansour, T.}
\newblock Restricted {M}otzkin permutations, {M}otzkin paths, continued
  fractions and {C}hebyshev polynomials.
\newblock {\em Discrete Math. 305}, 1-3 (2005), 170--189.

\bibitem{giakoumakis:on-the-closure-:}
{\sc Giakoumakis, V.}
\newblock On the closure of graphs under substitution.
\newblock {\em Discrete Math. 177}, 1-3 (1997), 83--97.

\bibitem{guibert:restricted-132-:}
{\sc Guibert, O., and Mansour, T.}
\newblock Restricted $132$-involutions.
\newblock {\em S\'em. Lothar. Combin. 48\/} (2002), Art.\ B48a, 23 pp.
  (electronic).

\bibitem{guibert:some-statistics:}
{\sc Guibert, O., and Mansour, T.}
\newblock Some statistics on restricted $132$ involutions.
\newblock {\em Ann. Comb. 6}, 3-4 (2002), 349--374.

\bibitem{gp:vex}
{\sc Guibert, O., and Pergola, E.}
\newblock Enumeration of vexillary involutions which are equal to their
  mirror/complement.
\newblock {\em Discrete Math. 224}, 1-3 (2000), 281--287.

\bibitem{higman:ordering-by-div:}
{\sc Higman, G.}
\newblock Ordering by divisibility in abstract algebras.
\newblock {\em Proc. London Math. Soc. (3) 2\/} (1952), 326--336.

\bibitem{knuth:the-art-of-comp:}
{\sc Knuth, D.~E.}
\newblock {\em The art of computer programming. {V}ol. 1: {F}undamental
  algorithms}.
\newblock Addison-Wesley Publishing Co., Reading, Mass., 1969.

\bibitem{lakshmibai:criterion-for-s:}
{\sc Lakshmibai, V., and Sandhya, B.}
\newblock Criterion for smoothness of {S}chubert varieties in {${\rm
  SL}(n)/B$}.
\newblock {\em Proc. Indian Acad. Sci. Math. Sci. 100}, 1 (1990), 45--52.

\bibitem{mansour:restricted-1-3-:}
{\sc Mansour, T.}
\newblock Restricted $1\mbox{\small-}3\mbox{\small-}2$ permutations and
  generalized patterns.
\newblock {\em Ann. Comb. 6}, 1 (2002), 65--76.

\bibitem{mansour:restricted-132-:}
{\sc Mansour, T.}
\newblock Restricted $132$-alternating permutations and {C}hebyshev
  polynomials.
\newblock {\em Ann. Comb. 7}, 2 (2003), 201--227.

\bibitem{mansour:counting-occurr:a}
{\sc Mansour, T.}
\newblock Counting occurrences of $132$ in an even permutation.
\newblock {\em Int. J. Math. Math. Sci.}, 25-28 (2004), 1329--1341.

\bibitem{mansour:restricted-132-:b}
{\sc Mansour, T.}
\newblock Restricted $132$-{D}umont permutations.
\newblock {\em Australas. J. Combin. 29\/} (2004), 103--117.

\bibitem{mansour:restricted-even:}
{\sc Mansour, T.}
\newblock Restricted even permutations and {C}hebyshev polynomials.
\newblock {\em Discrete Math. 306}, 12 (2006), 1161--1176.

\bibitem{mansour:restricted-132-:a}
{\sc Mansour, T., and Vainshtein, A.}
\newblock Restricted $132$-avoiding permutations.
\newblock {\em Adv. in Appl. Math. 26}, 3 (2001), 258--269.

\bibitem{mansour:counting-occurr:}
{\sc Mansour, T., and Vainshtein, A.}
\newblock Counting occurrences of $132$ in a permutation.
\newblock {\em Adv. in Appl. Math. 28}, 2 (2002), 185--195.

\bibitem{mansour:counting-occurr:b}
{\sc Mansour, T., Yan, S. H.~F., and Yang, L. L.~M.}
\newblock Counting occurrences of $231$ in an involution.
\newblock {\em Discrete Math. 306}, 6 (2006), 564--572.

\bibitem{mohring:algorithmic-asp:a}
{\sc M{\"o}hring, R.~H.}
\newblock Algorithmic aspects of the substitution decomposition in optimization
  over relations, sets systems and {B}oolean functions.
\newblock {\em Ann. Oper. Res. 4}, 1-4 (1985), 195--225.

\bibitem{mohring:substitution-de:}
{\sc M{\"o}hring, R.~H., and Radermacher, F.~J.}
\newblock Substitution decomposition for discrete structures and connections
  with combinatorial optimization.
\newblock In {\em Algebraic and combinatorial methods in operations research},
  vol.~95 of {\em North-Holland Math. Stud.} North-Holland, Amsterdam, 1984,
  pp.~257--355.

\bibitem{murphy:restricted-perm:}
{\sc Murphy, M.~M.}
\newblock {\em Restricted permutations, antichains, atomic classes, and stack
  sorting}.
\newblock PhD thesis, Univ. of St Andrews, 2002.

\bibitem{robertson:refined-restric:}
{\sc Robertson, A., Saracino, D., and Zeilberger, D.}
\newblock Refined restricted permutations.
\newblock {\em Ann. Comb. 6}, 3-4 (2002), 427--444.

\bibitem{schmerl:critically-inde:}
{\sc Schmerl, J.~H., and Trotter, W.~T.}
\newblock Critically indecomposable partially ordered sets, graphs, tournaments
  and other binary relational structures.
\newblock {\em Discrete Math. 113}, 1-3 (1993), 191--205.

\bibitem{shapiro:bootstrap-perco:}
{\sc Shapiro, L., and Stephens, A.~B.}
\newblock Bootstrap percolation, the {S}chr{\"o}der numbers, and the $n$-kings
  problem.
\newblock {\em SIAM J. Discrete Math. 4}, 2 (1991), 275--280.

\bibitem{stanley:enumeration-of-:}
{\sc Stanley, R.~P.}
\newblock Enumeration of posets generated by disjoint unions and ordinal sums.
\newblock {\em Proc. Amer. Math. Soc. 45\/} (1974), 295--299.

\bibitem{stanley:enumerative-com:}
{\sc Stanley, R.~P.}
\newblock {\em Enumerative combinatorics. {V}ol. 1}, vol.~49 of {\em Cambridge
  Studies in Advanced Mathematics}.
\newblock Cambridge University Press, Cambridge, 1997.

\bibitem{stanley:enumerative-com:a}
{\sc Stanley, R.~P.}
\newblock {\em Enumerative combinatorics. {V}ol. 2}, vol.~62 of {\em Cambridge
  Studies in Advanced Mathematics}.
\newblock Cambridge University Press, Cambridge, 1999.

\bibitem{west:sorting-twice-t:}
{\sc West, J.}
\newblock Sorting twice through a stack.
\newblock {\em Theoret. Comput. Sci. 117}, 1-2 (1993), 303--313.

\bibitem{zverovich:a-finiteness-th:}
{\sc Zverovich, I.}
\newblock A finiteness theorem for primal extensions.
\newblock {\em Discrete Math. 296}, 1 (2005), 103--116.

\end{thebibliography}

\end{document}